\documentclass[12pt]{amsart}

\usepackage{amssymb,amsmath,bm, verbatim, natbib, epsf, here, amsfonts,latexsym, graphicx}

\graphicspath{{Plots-NullDists/}}
\newtheorem{theorem}{Theorem}

\newcommand{\UU}{{\mathbb U}}

\numberwithin{equation}{section}

\begin{document}


\title{An excursion approach to maxima of the Brownian Bridge}
\author{Mihael Perman}
\author{Jon A. Wellner}
\date{\today}
\keywords{Brownian bridge, rescaling, excursions, extrema, \\Kolmogorov-Smirnov statistics.\hfill}
\address{Mihael Perman,
Faculty of Mathematics and Physics,
Jadranska 19,
SI-1000 Ljubljana,
Slovenia and
Faculty of Mathematics, Natural Sciences and Information Technologies,
University of Primorska, Glagolja\v ska 8, SI-6000 Koper, Slovenia.}
\email{mihael.perman@fmf.uni-lj.si}
\address{Jon A. Wellner, Department of Statistics, University of Washington,
B320 Padelford Hall, Seattle, WA 98195-4322, U.S.A.}
\email{jaw@stat.washington.edu}

\begin{abstract}
Distributions of functionals of Brownian bridge arise as limiting distributions in
nonparametric statistics.  In this paper we will give a derivation of
distributions of extrema of the Brownian bridge based on excursion theory
for Brownian motion.  The idea of rescaling and conditioning on the local
time has been used widely in the literature.  In this paper it is used to
give a unified derivation of a number of known distributions, and a few new
ones.  Particular cases of calculations include the 
distribution of the Kolmogorov-Smirnov statistic and the Kuiper statistic.

\end{abstract}

\maketitle

\newpage
\section{Introduction.}

Distributions of functionals of Brownian bridge arise as limiting distributions
in non-parametric statistics.  The distribution of the maximum
of the absolute value of a Brownian bridge is the basis for the Kolmogorov-Smirnov
non-parametric test of goodness of fit to give one example.  For an overview of statistical 
applications see \cite{MR838963}. 

Let $(\UU_t: 0 \le t \le 1)$ be the standard Brownian bridge and define
\begin{equation}
\label{DefM+m-}
M^+ = \max_{0 \le t \le 1} \UU_t, \qquad 
M^- = -\min_{0 \le t \le 1} \UU_t 
\end{equation}
and
\begin{equation}
\label{DefmM}
m=\min\{M^+,M^-\}
\qquad
\mbox{and}
\qquad
M=\max\{M^+,M^-\} .
\end{equation}
The distribution of $M^+$ was first computed by \cite{MR0001483}.  
The derivation of the distribution of $M$ was given by \cite{kol33}. 
For an elementary derivation
see \cite{MR0348887}.  
In this paper we give a derivation of the joint distribution of $M^+$ and $M^-$ based
on exursion theory for Brownian motion.  The distributional results can be used to
derive known distributions like the distribution of the Kuiper statistic
$K=M^++M^-$, or the distribution of the difference $D=M^+-M^-$ which seems to be new.

Let $(B_t \colon t \ge 0)$ be standard Brownian motion.
Define the last exit time from $0$ of $B$ before time $t$ as
\begin{equation}
\label{LastExitTime}
g_t = \sup\{ s \le t \colon B_s = 0 \} \,.
\end{equation}
The following lemma is well known and will be used to derive distributional equalities
needed later.
See \cite{MR0029120},  
\cite{MR0116376},  
\cite{MR1022918}.  

\begin{theorem}
\label{Dynkin}
The distribution of $g_1$ is
$\mathrm{Beta}(1/2,1/2)$.  Given $g_1$, the process $(B_t:0 \le t \le g_1)$
is a Brownian bridge of length $g_1$, and the randomly rescaled process
$(B(g_1 u)/\sqrt{g_1}: 0\le u\le 1)$ is a Brownian bridge independent
of $g_1$.
\end{theorem}
Let $S_\theta \sim \exp(\theta)$ be independent from $B$.  From scaling properties
of Brownian motion and Theorem 1 it follows that the process
\begin{equation}
\label{BMuptoS}
\left ( \frac{B_{tg_{S_\theta}}}{\sqrt{g_{S_\theta}}} \colon 0 \le t \le 1 \right )
\end{equation}
is a Brownian bridge independent of $g_{S_\theta}$.  Furthermore, the law of
$g_{S_\theta}$ is equal to the law of $g_1 S_\theta$ where $g_1$ and $S_\theta$
are assumed to be independent which is known to be $\Gamma(1/2,\theta)$; 
here $X \sim \Gamma (a,b)$ means that $X$ has the Gamma density 
\begin{equation}
p(x; a,b) = \frac{b (bx)^{a-1}}{\Gamma (a)} \exp (-bx) 1_{(0,\infty)} (x) .\nonumber
\end{equation}

Let $\UU$ be the standard Brownian bridge.
Let $\gamma \sim \Gamma(1/2,\theta)$ be independent from $\UU$ and let
\begin{equation}
\label{RescaledBridge}
\tilde \UU_t = \sqrt{\gamma} \, \UU_{t/\gamma}
\end{equation}
for $0 \le t \le \gamma$.  The process $(\tilde \UU_t \colon 0 \le t \le \gamma)$ is
called the randomly rescaled Brownian bridge.  From the independence of $g_{S_\theta}$ and
the process defined in (\ref{BMuptoS}) we have 
\begin{equation}
\label{EqualityinLaw}
(\tilde \UU_t \colon 0 \le t \le \gamma) \stackrel{\mathrm{d}}{=}
\left (B_t \colon 0 \le t \le g_{S_\theta} \right ) \,.
\end{equation}
This equality in law can be exploited to derive Laplace transforms of distributions 
of functionals of Brownian bridge.  From (\ref{EqualityinLaw}) it follows that
\begin{equation}
\label{BasicEqualityinLaw}
\left (\sqrt{\gamma} M^+,\sqrt{\gamma} M^- \right )
\stackrel{\mathrm{d}}{=} 
\left ( \max_{0 \le t \le g_{S_\theta}} B_t, -\min_{0 \le t \le g_{S_\theta}} B_t \right ) \,.
\end{equation}
Excursion theory will provide the distribution of the pair on the  right in
(\ref{BasicEqualityinLaw}) which in turn is used to derive the Laplace transform of the
cumulative distribution function of the pair $(M^+,M^-)$.  The transform can be inverted
in the form of infinite series.
The method has been used widely in the literature and is well known.  
See \cite{py97kt}, \cite{cppy98}, \cite{py98csaki} and \cite{MR1422979} for
results based on this identity in law.  
The contribution of this paper is a unified way to derive explicitly the distributions 
of functionals related to the pair $(M^+,M^-)$.

\section{Brownian excursions.}

The paths of Brownian motion $B$ are continuous functions hence the complements
of their zero sets are unions of disjoint open intervals.  The path of Brownian motion
restricted to any such open interval is called an excursion away from $0$.
Since Brownian motion is recurrent all open intervals will be bounded.
The path can thus be broken up into an infinite string of excursions and
every excursion can be identified with a function in the set of functions
\begin{equation}
\label{ExcursionSpace}
\mathcal{U}=\{ w \in C[0,\infty), w(0)=0, \exists R > 0, w(t) \ne 0 \;{\rm iff} \; t\in (0,R) \} \,.
\end{equation}
To describe the structure of excursions let $(L(t):t \ge 0)$ be the local time 
process at level 0 for Brownian motion normalized so that
\begin{equation}
\label{LocalTimeStandardization}
L(t) \stackrel{\mathrm{d}}{=} \max_{0 \le s \le t} \, B_s \,.
\end{equation}
Local time is an adapted nondecreasing process which only increases on the zero set
of Brownian motion and $L(t_1) < L(t_2)$ for $t_1 < t_2$ whenever the interval
$(t_1,t_2)$ contains a zero of $B$.  Hence the local time during two different excursions
is different and constant during each excursion.  See \cite{MR1725357}  
for definitions and fundamental results on local time.
Let $\tau_s = \inf\{u:L(u) > s\}$ be the right continuous inverse of
the local time process $L(t)$.  From the properties of local time we infer that
every excursion of Brownian motion away from $0$ is of the form
\begin{equation}
\label{excursions}
e_s(u) = 1_{[0 \le u \le  \tau_s-\tau_{s-})}(u) \, B_{\tau_{s-} + u} 
\end{equation}
for those $s$ at which $\tau_s$ has a jump.  Let $e$ be the point process
defined on the abstract space $(0,\infty) \times  \mathcal{U}$ defined as
\begin{equation}
\label{DefinitionOfExcursions}
e=\{(s,e_s) \colon s>0, \tau_s-\tau_{s-} > 0\} \,.
\end{equation}
The following theorem by It\^o is one of his great insights.
\begin{theorem}
\label{ItosExcursions}
The point process $e$ is a Poisson process on $(0,\infty) \times \mathcal{U}$ with mean measure
given by $\lambda \times n$ where $\lambda$ is the Lebesgue measure on $(0,\infty)$
and $n$ is a $\sigma$-finite measure on the functions space $\mathcal{U}$ equipped with the
$\sigma$-field generated by the coordinate maps.
\end{theorem}
\noindent
For a proof see  \cite{MR1725357}, p. 457.  
Note that the excursions of the process $(B_t \colon 0 \le t \le g_{S_\theta})$ are
a portion of the excursion process of Brownian motion.   It will be shown that the law
of this portion can be described and used to derive the distribution of the pair
of variables on the right side of (\ref{BasicEqualityinLaw}).  

If the points of a Poisson process on an abstract space with mean
measure $\mu$ are
marked in such a way that each point receives a mark independently of other points
with probability depending on the position of the point then the point processes
of marked and unmarked excursions are two independent Poisson processes.  If at position
$x$ a mark is assigned with probability $f(x)$ then the marked and unmarked Poisson processes
have mean masures $f \cdot \mu$ and $(1-f)\cdot \mu$.  See
\cite{MR1207584}, p. 55, for definitions and proof.  

Marking will be applied to the Poisson process of excursions.  Define the duration of an 
excursion $w \in \mathcal{U}$ as
\[
R(w)=\sup\{u \colon w(u) \ne 0\} \,.
\]
and assign marks to the points of the process of excursions with probability
$1-e^{-\theta R(w)}$ for $\theta>0$.  Define $T=\inf\{s \colon \mbox{$e_s$ is marked}\}$.  

\begin{theorem}
\label{markedexcursions}
Let $\tilde e$ be the point process $\{(s,e_s) \colon 0 < s <T\}$.

\begin{itemize}

\item[(i)] The random variable $T$ is exponential with parameter $\sqrt{2\theta}$.
\item[(ii)] Conditionally on $T=t$, the point process $ \tilde e$ is a Poisson process
in the space  $(0,t) \times U$ with mean measure 
$\lambda \times n\cdot e^{-\theta R(w)}$
where $\lambda$ is the Lebesgue measure on $(0,t)$ and $n$ is It\^o's excursion law for 
Brownian motion.
\item[(iii)] Positive and negative excursions of $\tilde e$ are conditionally 
independent Poisson processes given  $T=t$.

\end{itemize}

\end{theorem}

\noindent
\begin{proof}[{\sc Proof}:]  For the It\^o measure $n$ we have
\[
n(R \in \mathrm{d}r)=\frac{\mathrm{d}r}{\sqrt{2\pi r^3}} \,.
\]
See \cite{MR1725357}, p. 459.  The probability that there is no marked excursion 
in $(0,t) \times \mathcal{U}$ is given by
\[
\exp \left (-\int_{(0,t) \times U} \, (1-e^{-\theta r}) \, \mathrm{d}t \, n( R \in \mathrm{d}r) \right )
=
e^{-t \sqrt{2\theta} } \,.
\]
The integral above is given e.g. in \cite{MR1406564}, p. 73.
It follows that $P(T > t) = e^{-t \sqrt{2\theta} }$ which proves (i).  The assertions
in (ii) and (iii) follow from independence properties of Poisson processes bearing in mind
that positive and negative excursions of $e$ are independent Poisson processes and that
marking is independent of the sign of excursions.
\end{proof}

The law of excursions of $(B_t \colon 0 \le t \le g_{S_\theta})$ is described in the following 
theorem.
See also \cite{MR921238}, p. 418.

\begin{theorem}
\label{markedexcursions2}
The law of the point process $((s,e_s) \colon 0 < s < L(S_\theta))$ is described by:

\begin{itemize}
\item[(i)] $L(S_\theta)$ is exponential with parameter $\sqrt{2\theta}$.
\item[(ii)] Conditionally on $L(S_\theta)=t$ the process $((s,e_s) \colon 0 < s < L(S_\theta))$
is a Poisson process on $(0,t) \times \mathcal{U}$ with mean measure
$\lambda \times n \cdot e^{-\theta R(w)}$.
\item[(iii)] Conditionally on $L(S_\theta)=t$ the positive and negative excursions of
$((s,e_s) \colon 0 < s < L(S_\theta))$ are independent Poisson processes.
\end{itemize}
\end{theorem}

\noindent
\begin{proof}[{\sc Proof}:] Let $N$ be a Poisson process with intensity $\theta$ on $(0,\infty)$ independent
of $B$.  If $e_s$ is an excursion on the open interval of length $R$ then by independence
the open interval contains a point of $N$ with probability $1-e^{-\theta R}$.  Declare all
excursions that contain a point of $N$ to be marked.  By independence properties of $N$ marks 
are assigned independently.  The leftmost point of $N$ will be an exponential random variable
$S_\theta$ independent of $B$.  It follows that the excursion straddling $S_\theta$ is
exactly the first marked excursion of the excursion process $e$.  Hence the excursion process of
$(B_t \colon 0 \le t \le g_{S_\theta})$ is exactly the portion of the excursion process $e$
up to the first marked excursion.  The assertions (i), (ii) and (iii) follow 
from Theorem \ref{markedexcursions}.
\end{proof}

Let $M^+$, $M^-$ and $M$ be as defined in Section 1.  Let $\UU$ be
a Brownian bridge and $\gamma$ a $\Gamma(1/2,\theta)$ random
variable independent of $\UU$.  Some preliminary calculations are needed to find
explicitly the distribution of $(\sqrt{\gamma} M_+,\sqrt{\gamma} M^-)$.

The reflection principle for Brownian motion states, see
\cite{MR1912205}, p. 126, formula 1.1.8, that
for $x > 0$ and $z < y < x$
\begin{equation}
\label{ReflectionPrinciple}
P \left ( \max_{0\le t\le 1} \, B_t \ge x, z < B_1 < y \right )
=
P \left ( 2x-y < B_1 < 2x-z \right ) \,.
\end{equation}
Brownian bridge is Brownian motion conditioned
to be 0 at time $t=1$ so by (\ref{ReflectionPrinciple}) for $x>0$
\begin{eqnarray}
\label{onesided}
P( M^+ \ge x ) 
& = &
\lim_{\epsilon\rightarrow 0} 
P \left (\max_{0\le t\le 1} \, B_t \ge x \, \big \vert \, \vert B_1 \vert \le \epsilon \right )
\nonumber \\
& = &
\lim_{\epsilon\rightarrow 0} 
P \left (\max_{0\le t\le 1} B_t \ge x , \vert B_1\vert \le \epsilon \right )/
P \left (\vert B_1\vert \le \epsilon  \right ) \nonumber \\
& = &
\lim_{\epsilon\rightarrow 0} 
P \left (2x-\epsilon \le B_1 \le 2x+\epsilon  \right )/
P \left (\vert B_1\vert \le \epsilon  \right ) \nonumber \\
& = &
e^{-2x^2} \,. \\ \nonumber
\end{eqnarray}
It follows from the distribution of $M^+$ given by (\ref{onesided}) that
\begin{eqnarray}
\label{rescaledm+}
P(\sqrt{\gamma} M^+ \ge x) 
& = & \sqrt{\theta/\pi} \int_0^\infty  \exp(\,-2x^2/s ) s^{-1/2} e^{-\theta s} {\rm d}s \\
& = & \exp(\,-2x\sqrt{2\theta}) \,. \nonumber 
\end{eqnarray}
The integral is given in \cite{MR0352889}, p. 41,
formula 5.28.

Turning to excursions recall that $R(w)$ stands for the
length of the excursion and denote $w^+=\max_u w(u)$.  Define for $x>0$
\[
m(x) = \int_{\{w^+ > x\}} e^{-\theta R(w)}  n(\mathrm{d}w) \,.
\]
\begin{theorem}
\label{Triple}
The law of the triple 
\[
\left (  \max_{0 \le t \le g_{S_\theta}} B_t, -\min_{0 \le t \le g_{S_\theta}} B_t,L(S_\theta) \right )
\]
is described by:
\begin{itemize}
\item[(i)] $L(S_\theta)$ is exponential with parameter $\sqrt{2\theta}$.
\item[(ii)] The random variables $\max_{0 \le t \le g_{S_\theta}} B_t$ and
$ -\min_{0 \le t \le g_{S_\theta}} B_t$ are conditionally independent given
$L(S_\theta)=t$ with the same conditional distribution.
\item[(iii)] 
\[
P \left ( \max_{0 \le t \le g_{S_\theta}} B_t \le  x \vert L(S_\theta)=t \right )
=
e^{-t m(x)} \,.
\]
\item[(iv)] 
\begin{equation}
\label{Functionm}
m(x) = \frac{\sqrt{2\theta} e^{-2x\sqrt{2\theta}}}{1-e^{-2x\sqrt{2\theta}}} \,.
\end{equation}
\end{itemize}
\end{theorem}

\noindent
\begin{proof}[{\sc Proof}:] (i) is proved in Theorem \ref{markedexcursions2}.
The point processes of positive and negative excursions of
$((s,e_s) \colon 0 < s < L(S_\theta))$ are conditionally independent
given $L(S_\theta)=t$ by Theorem \ref{markedexcursions2}.  Since 
$\max_{0 \le t \le g_{S_\theta}} B_t$ is a function of positive and
$ -\min_{0 \le t \le g_{S_\theta}} B_t$ a function of negative excursions
conditional independence follows.  Equality of conditional distributions follows by
symmetry.  The process
$\left ((s, e_s^+) \colon 0 <s < L(S_\theta) \right )$ 
is a measurable map of the process $\left ((s,e_s) \colon 0 < s < L(S_\theta) \right )$
hence conditionally on $L(S_\theta)=t$ a Poisson process on $(0,t)\times (0,\infty)$;  see 
\cite{MR1207584}.  Conditionally on $L(S_\theta)=t$ we have that
$\max_{0 \le t \le g_{S_\theta}} B_t \le  x$ if there is no point of
$\left ((s, e_s^+) \colon 0 <s < L(S_\theta) \right )$ in the set $(0,t)\times (x,\infty)$.
The measure of this set is $t m(x)$ by Theorem \ref{markedexcursions2}, (ii).  The assertion
(iii) follows.  By unconditioning
\begin{eqnarray}
P( \max_{0 \le t \le g_{S_\theta}} B_t > x )  
& = &  \sqrt{2\theta} \int_0^\infty e^{-\sqrt{2\theta} t} (1-e^{-t m(x)}) \, \mathrm{d}t \label{excursionm+} \\
& = &  \frac{m(x)}{\sqrt{2\theta} + m(x)} \,.  \nonumber
\end{eqnarray}
Comparing (\ref{rescaledm+}) and (\ref{excursionm+}) we obtain that
\begin{equation}
\frac{m(x)}{\sqrt{2\theta} + m(x)} = \exp \left (\,-2x\sqrt{2\theta} \right )\,.
\end{equation}
(iv) follows by solving for $m(x)$.
\end{proof}
The law of the triple is in accordance with formula (53) in \cite{py98csaki}.

\section{Examples of calculations}

\leftline{\sc 3.1 Distributions of $M$ and $m$.}

\medskip
Let $m$ and $M$ be defined as in (\ref{DefmM}).  By (\ref{BasicEqualityinLaw}) the random variables
$\max_{0 \le t \le g_{S_\theta}} B_t$ and $\sqrt{\gamma} M^+$ have the same distribution.
We compute, by conditioning on $L(S_{\theta})$ and using Theorem~\ref{Triple},
\begin{eqnarray}
P( \sqrt{\gamma} M \le x ) 
& = & P( \sqrt{\gamma} M^+ \le x , \ \sqrt{\gamma} M^{-} \le x) \nonumber \\
& = & P( \max_{0 \le t \le g_{S_\theta}} B_t \le x, \ - \min_{0 \le t \le g_{S_\theta}} B_t \le x) \nonumber \\
& = & E \{ P( \max_{0 \le t \le g_{S_\theta}} B_t \le x, \ - \min_{0 \le t \le g_{S_\theta}} B_t \le x | L(S_{\theta} ) ) \} \nonumber \\
& = & E \{ \exp ( - 2 L(S_{\theta} ) m(x) ) \} \nonumber \\
& = & \int_0^{\infty} \sqrt{2 \theta} \exp (- \sqrt{2\theta} v) \exp ( - 2 v m(x) ) dx \nonumber  \\
& = & \frac{\sqrt{2\theta}}{\sqrt{2\theta} + 2 m(x)}  = \tanh (x \sqrt{2 \theta} )  \label{LaplaceMDeriv}
\end{eqnarray}
where we used (\ref{Functionm}) in the last step.
Thus 
\begin{equation}
\label{LaplaceM}
P(\sqrt{\gamma} M \le x) =
\frac{\sqrt{2\theta}}{\sqrt{2\theta} + 2 m(x)} = \tanh(x\sqrt{2\theta}) \,.
\end{equation}
Let $F_M$ be the cumulative distribution function of $M$.  Writing out (\ref{LaplaceM}),
taking into account the independence of $M$ and $\gamma$ and
dividing both sides by $\sqrt{\theta}$ we get
\begin{equation}
\label{oberhettinger8.51}
\frac{1}{\sqrt{\pi}} \,  \int_0^\infty F_M(x/\sqrt{s}) s^{-1/2} e^{-\theta s} {\rm d}s =
\frac{\tanh(x\sqrt{2\theta})}{\sqrt{\theta}} \,.
\end{equation}
\cite{MR0352889}, p. 294, 
formula 8.51, 
give the inverse of the Laplace transform on the right of (\ref{oberhettinger8.51}) as
\begin{eqnarray*}
\frac{1}{\sqrt{\pi}} F_M(x/\sqrt{s}) s^{-1/2}
& = &
\frac{1}{\sqrt{2} x} \theta_2(0\vert s/(2x^2)) \\
& = &
\frac{1}{\sqrt{2} x} \frac{\sqrt{2}x}{\sqrt{\pi s}}
\sum_{k=-\infty}^\infty (-1)^k \exp(\,- 2 k^2 x^2/s ) \,. \\
\end{eqnarray*}
This yields
\[
F_M(z) = 
\sum_{k=-\infty}^{\infty} (-1)^k \exp(\,- 2 k^2 z^2 ),
\]
or
\[
1-F_M(z) = 2 
\sum_{k=1}^{\infty} (-1)^{k+1} \exp(\,- 2 k^2 z^2 ), 
\]
which is the formula for the distribution of the Kolmogorov--Smirnov test
statistic.

Turning to $m=\min\{M^+,M^-\}$ observe that by Theorem \ref{Triple}
\begin{equation}
\label{minimumM=M-}
P \left ( \sqrt{\gamma} \, m > x \right )
=
\sqrt{2\theta} \int_0^\infty \, e^{-\sqrt{2\theta} t} \, \left ( 1-e^{-tm(x)} \right )^2 \, \mathrm{d}t \,.
\end{equation}
Integration yields
\begin{equation}
\label{Laplacem}
P \left ( \sqrt{\gamma} \, m > x \right )
=
1-\frac{2\sqrt{2\theta}}{\sqrt{2\theta}+m(x)} + \frac{\sqrt{2\theta}}{\sqrt{2\theta}+2m(x)} \,.
\end{equation}
Denote by $F_m$ the distribution function of $m$.  From (\ref{Laplacem}), independence of
$m$ and $\gamma$ and dividing both sides by $\sqrt{\theta}$ we get
\[
\frac{1}{\sqrt{\pi}} \int_0^\infty \left (1-F_m(x/\sqrt{s}) \right ) s^{-1/2} e^{-\theta s} {\rm d}s =
\frac{\tanh(x\sqrt{2\theta})}{\sqrt{\theta}} +
\frac{2e^{-2x\sqrt{2\theta}}}{\sqrt{\theta}} - \frac{1}{\sqrt{\theta}} \,.
\]
The first term has been inverted above, the second is given by \cite{MR0352889}, p. 258, formula
5.87, and the third is elementary.  Substituting $z$ for $x/\sqrt{s}$ one gets
\begin{eqnarray}
\label{Fm}
1-F_m(z)
& = &
\sum_{k=-\infty}^\infty (-1)^k \exp(\,- 2 k^2 z^2 ) +
2e^{-2z^2} - 1 \nonumber \\
& = &
2 \, \sum_{k=2}^\infty \,  (-1)^k \exp(\,- 2 k^2 z^2 ) \,.  \\ \nonumber
\end{eqnarray}

\leftline{\sc 3.2 Joint distributions, sums, differences, quotients.}

\medskip
In this section the distributions of various functions of the pair $(M^+,M^-)$ will
be derived.  Let $K\equiv M^{+} + M^{-}$, the Kuiper (or range) statistic, 
$L \equiv M^{+} - M^{-}$, the difference statistic, and let 
$Q \equiv M^{+} / M^{-} $, the ratio statistic.  

\begin{theorem}  
\label{SumsDifferencesQuotients}
$\phantom{blab}$\\
(i) \ \ 
The joint distribution of $(M^+, M^-)$ is given for $x,y>0$ by
\begin{eqnarray}
\lefteqn{P(\sqrt{\gamma} M^+ \le x,\sqrt{\gamma} M^- \le y) } \nonumber  \\
& = & \frac{\sqrt{2\theta}}{\sqrt{2\theta} +m(x) + m(y)} \, \nonumber \\
& = & \coth((x+y)\sqrt{2 \theta}) 
      - \frac{\cosh ((x-y)\sqrt{2 \theta})}{ \sinh ((x+y)\sqrt{2 \theta})} \,.
      \label{JointRescaledStatement} 
\end{eqnarray}
(ii) \ \ For $x>0$,
\begin{equation}
\label{Kuiper}
F_K(x)  = P(M^{+} + M^{-} \le x) = 
\sum_{k=-\infty}^\infty (1 - 4k^2 x^2 ) e^{-2k^2x^2} \,.
\end{equation}
(iii) \ For $x>0$,
\begin{equation}
1-F_L(x) = P(M^+ - M^- \ge x) = \sum_{k=1}^\infty \frac{1}{4k^2-1} \cdot e^{-2k^2 x^2} \,.
\end{equation}
(iv) \ For $x>0$,
\begin{equation}
F_Q(x) = P( Q \le x) = \frac{1}{z+1} \, \left ( 1-\frac{\pi z \mathrm{cot}\left ( \frac{\pi z}{z+1} \right )}{z+1} \right ) .
\end{equation}
\end{theorem}

\medskip
\noindent
{\sc Remark 1:}  The formula for the joint distribution  in (i) is
in agreement with 
\cite{MR838963},  
formula (2.2.22), page 39.
\medskip

\noindent
{\sc Remark 2:}  
The result (ii) is in agreement with \cite{MR0111088} and with 
\cite{MR0488202}, Proposition 22.10, page 22.6.
\cite{MR515820} 
gives a construction of standard Brownian
excursion from a Brownian bridge.  Let $\UU$ be a Brownian bridge $[0,1]$ and
let $\sigma$ be the time when $\UU$ attains its minimum on $[0,1]$ ($\sigma$ is
a.s. unique).  Then the process $(e(t) \colon 0 \le t \le 1)$ defined by
\[
e(t) = \UU_{\sigma+t\; {(mod\; 1)}} - \UU_{\sigma}(t) 
\]
is a standard Brownian excursion.  It is a simple consequence of this
transformation that the Kuiper statistic has the distribution of the
maximum of the standard Brownian excursion and (\ref{Kuiper}) is another
derivation of the distribution of this maximum.  
Further results for range statistics are given by \cite{MR0042626} and \cite{MR-SV:07}.
\medskip

\par\noindent
{\sc Remark 3:}  
Note that the distribution of $M^{+} - M^{-}$ is symmetric about $0$.
For a different approach for Brownian motion instead of Brownian bridge see \cite{MR1689236}.  
\medskip

\par\noindent
{\sc Remark 4:} The result in (iv) is in accordance
with the distribution for the ratio $\widetilde{Q}=M^+/(M^++M^-)$ given in \cite{MR561878}.  The
derivation of the distribution of $\widetilde{Q}$ based on rescaling arguments is given in 
\cite{py98csaki} and can be derived easily from the above result for $F_Q$. 

\noindent
\begin{proof}[{\sc Proof}:] 
By unconditioning in Theorem \ref{Triple}
\begin{eqnarray}
\label{JointRescaled}
\lefteqn{\quad P(\sqrt{\gamma} M^+ \le x,\sqrt{\gamma} M^- \le y) =}  \\
& = & \frac{\sqrt{2\theta}}{\sqrt{2\theta} +m(x) + m(y)} \, \nonumber \\
& = & \frac{(1 - e^{-2x\sqrt{2 \theta}})(1 - e^{-2y \sqrt{2 \theta}})}
       {1 - e^{-2x\sqrt{2\theta}} e^{-2y\sqrt{2\theta}} }  \nonumber  \\
& = & \frac{(e^{x\sqrt{2 \theta}} - e^{-x\sqrt{2 \theta}})
            ( e^{y\sqrt{2 \theta}} - e^{-y\sqrt{2 \theta}})}
           {e^{(x+y) \sqrt{2 \theta}} - e^{-(x+y) \sqrt{2 \theta}}} \nonumber \\
& = & 2 \frac{ \sinh(x \sqrt{2\theta}) \sinh(y \sqrt{2 \theta})}
             { \sinh((x+y) \sqrt{2 \theta})}  \nonumber \\
& = & \frac{\cosh((x+y)\sqrt{2 \theta}) - \cosh((x-y)\sqrt{2 \theta})}
           {\sinh((x+y)\sqrt{2 \theta})} \nonumber \\
& = & \coth((x+y)\sqrt{2 \theta}) 
      - \frac{\cosh ((x-y)\sqrt{2 \theta})}{ \sinh ((x+y)\sqrt{2 \theta})} \,.\nonumber \\ \nonumber
\end{eqnarray}      
Let $F(z,w)$ be the joint distribution function of the pair $(M^+,M^-)$.  By independence
of $\gamma$ and $(M^+,M^-)$
\begin{eqnarray}
\label{JointLaplaceTransform}
\lefteqn{P(\sqrt{\gamma} M^+ \le x,\sqrt{\gamma} M^- \le y)=}  \nonumber \\
& = &
\sqrt{\theta/\pi} \int_0^\infty F(x/\sqrt{s},y/\sqrt{s}) s^{-1/2} e^{-\theta s} {\rm d}s  \,. 
\end{eqnarray}
The right side is given in (\ref{JointRescaled}).  \cite{MR0352889}, p. 294, formula 8.52, 
give the inverse of the first term on the right in (\ref{JointRescaled})
\[
\frac{1}{\sqrt{2} (x+y)} \theta_3 \left (0\big \vert \frac{s}{2(x+y)^2} \right ) \,
= 
\frac{1}{\sqrt{\pi s}} 
\sum_{k=-\infty}^{\infty} e^{-2k^2 (x+y)^2 /s} 
\]
The inverse of the second term of the transform can be obtained from 
\cite{MR0352889}, p. 294, 
formula 8.60:  we find that the inverse is
\begin{eqnarray*}
\lefteqn{\frac{1}{\sqrt{2} (x+y)} \theta_4 \left ( \frac{(x-y)/2}{x+y} \big \vert \frac{s}{2(x+y)^2} \right )=} \\
& = & 
\frac{1}{\sqrt{\pi s}} 
\sum_{k=-\infty}^{\infty} 
\exp\left ( -2 (x+y)^2  \left [\frac{(x-y)/2}{(x+y)} + k + \frac{1}{2}  \right ]^2 /s \right )   \\
& = & 
\frac{1}{\sqrt{\pi s}} 
\sum_{k=-\infty}^{\infty} 
\exp( -2 [k(x+y) +x]^2/s ) \, .
\end{eqnarray*}
Combining these yields 
\begin{eqnarray*}
\lefteqn{P( M^+ \le z, M^- \le w) = } \\
& = &
\sum_{k=-\infty}^{\infty} \exp (-2k^2 (z+w)^2 )
- 
\sum_{k=-\infty}^{\infty} 
\exp \left ( -2 (k(z+w) +z)^2  \right ) \, .
\end{eqnarray*}

We now consider the Kuiper statistic  $K=M^++M^-$.
It seems cumbersome to proceed from the joint
distribution of $M^+$ and $M^-$ so we use directly the distribution
of $(\sqrt{\gamma}M^+,\sqrt{\gamma}M^-)$.  Denote $U=\sqrt{\gamma}M^+$ and
$V=\sqrt{\gamma}M^-$.  The joint cumulative distribution
function of $U$ and $V$ is given in (\ref{JointRescaled}) as

\[
G(u,v)= 2 \frac{\sinh(u\sqrt{2\theta})\sinh(v\sqrt{2\theta})}
               {\sinh((u+v)\sqrt{2\theta})}
\]
The cumulative distribution function of $U+V$ is given by

\[
P(U+V \le z) = \int_0^z G_u (u,z-u) {\rm d}u
\]
where $G_u$ is the partial derivative of $G$ with respect to $u$.  A
calculation yields

\[
G_u(u,z-u) = \frac{2\sqrt{2\theta}\sinh^2((z-u)\sqrt{2\theta})}
                  {\sinh^2(z\sqrt{2\theta})}
\]
and integration gives

\[
P(U+V \le z) = \coth(z\sqrt{2\theta}) -
               \frac{z\sqrt{2\theta}}{\sinh^2(z\sqrt{2\theta})} \,.
\]
Using the fact that $\gamma$ and $M^++M^-$ are independent we obtain
the Laplace transform of the cumulative distribution function $F_K$ of 
the Kuiper statistic as

\[
\frac{\sqrt{\theta}}{\sqrt{\pi}}
\int_0^\infty F_K(z/\sqrt{s}) s^{-1/2} e^{-\theta s} {\rm d}s =
               \coth(z\sqrt{2\theta}) -
               \frac{z\sqrt{2\theta}}{\sinh^2(z\sqrt{2\theta})} \,.
\]
After dividing by $\sqrt{\theta}$ it remains to invert the two terms on the
right and substitute for $z/\sqrt{s}$.  The first term has been inverted above
when deriving the joint distribution of $M^+$ and $M^-$.  We get

\[
\frac{1}{\sqrt{2} z} \theta_3\left (0\big \vert \frac{s}{2z^2} \right ) \,
= 
\frac{1}{\sqrt{\pi s}} 
\sum_{k=-\infty}^{\infty} e^{-2k^2 z^2 /s} \,.
\]
To invert the second term rewrite it as 

\begin{eqnarray*}
\frac{\sqrt{2} z}{\sinh^2(z\sqrt{2\theta})} \,
 = 
\frac{4\sqrt{2}z e^{-2z\sqrt{2\theta}}}{(1-e^{-2z\sqrt{2\theta}})^2}
 = 
4 \sqrt{2}z\,\sum_{k=1}^\infty k e^{-2kz\sqrt{2\theta}}
\end{eqnarray*}
for $z > 0$ and $\theta > 0$.  The inverse Laplace transforms of the terms 
in the sum are known , see \cite{MR0352889}, p. 258, formula 5.85.  Taking the derivative with respect 
to $x$ on both sides of (\ref{rescaledm+}) we get
\begin{equation}
\label{invert}
\int_0^\infty \frac{a}{\sqrt{2\pi s^3}} 
\exp(\,-a^2/2s) e^{-\theta s} {\rm d}s \,=
e^{-a\sqrt{2\theta}} \,.
\end{equation}
Since all the terms are nonnegative functions the order of summation and
integration can be changed.  Hence for $z > 0$ the inverse Laplace 
transform of the second term is

\[
4 \sqrt{2} z \sum_{k=1}^\infty 
\frac{2k^2 z}{\sqrt{2 \pi s^3}} \exp(\,-2k^2 z^2/s) 
= \frac{1}{\sqrt{\pi s}} \sum_{k=-\infty}^\infty \frac{4 k^2 z^2}{s} 
\exp(\,-2k^2 z^2/s) \,.
\]
Substitute $x=z/\sqrt{s}$ to get

\begin{equation}
\label{Kuiper1}
F_K(x) = \sum_{k=-\infty}^\infty e^{-2k^2x^2}
           -   \sum_{k=-\infty}^\infty 4 k^2x^2 e^{-2k^2x^2} 
= 
\sum_{k=-\infty}^\infty (1 - 4k^2 x^2 ) e^{-2k^2x^2}
\,.
\end{equation}

\medskip

For the difference $U-V$ a computation yields
for $z > 0$

\begin{equation}
\label{distributiondiff}
P( U - V \ge z ) = \int_z^\infty G_u(u,u-z) \,{\rm d}u 
\end{equation}
provided $P(U > 0) =1$ and $P(V > 0)=1$ which is the case for the
variables in question.  Using the joint cumulative distribution
function yields
\[
G_u(u,u-z) = 2 \sqrt{2\theta}
\frac{\sinh^2(\sqrt{2\theta}(u-z))}{\sinh^2(\sqrt{2\theta}(2u-z))} \,.
\]
From(\ref{distributiondiff}) it follows
\begin{eqnarray*}
\label{DiffUV}
P( U - V \ge z )  
& = &
2 \sqrt{2\theta} \int_z^\infty 
\frac{\sinh^2(\sqrt{2\theta}(u-z))}{\sinh^2(\sqrt{2\theta}(2u-z))} \,{\rm d}u \\
& = &
2 \sqrt{2\theta} \int_0^\infty 
\frac{\sinh^2(\sqrt{2\theta}u)}{\sinh^2(\sqrt{2\theta}(2u+z))} \,{\rm d}u \,.\\ 
& = &
\frac{1}{2} - \mathrm{arctanh}(e^{-\sqrt{2\theta} z}) \sinh(\sqrt{2\theta} z) \,. 
\end{eqnarray*}
The integral in the last line is elementary and is computed by Mathematica.
Since $\gamma$ and $(M^+,M^-)$ are independent we have for $z>0$
\begin{eqnarray}
\label{LTOfU-V}
\lefteqn{P(U-V \ge z)=} \\ \notag
& = &
\frac{\sqrt{\theta}}{\sqrt{\pi}}
\int_0^\infty P \left (M^+ - M^- \ge \frac{z}{\sqrt{s}} \right )  \, 
\frac{e^{-\theta s} }{\sqrt{s}} \, \mathrm{d}s \\ \notag
& = &
\frac{1}{2} - \mathrm{arctanh}\left (e^{-\sqrt{2\theta} z} \right ) 
\sinh \left (\sqrt{2\theta} z \right ) \,. \\ \notag
\end{eqnarray}
This gives the Laplace transform of
$P(M^+ - M^- \ge z/\sqrt{s})/\sqrt{s}$ as a function of $s$ for fixed $z$.  

To invert this Laplace transform we use the known series expansion for the 
hyperbolic arc-tangent to get for $z > 0$
\begin{equation}
\label{ArcSeries}
\frac{1}{2} - \mathrm{arctanh}\left (e^{-\sqrt{2\theta} z} \right ) 
\sinh \left (\sqrt{2\theta} z \right )
=
\sum_{k=1}^\infty \, \frac{e^{-2kz\sqrt{2\theta}}}{4k^2-1} \,.
\end{equation}
Finally
\[
\int_0^\infty \, P \left (M^+-M^- \ge \frac{z}{\sqrt{s}} \right ) \, 
\frac{e^{-\theta s}}{\sqrt{s}} \mathrm{d}s
=
\sum_{k=1}^\infty \, \frac{(\pi/\theta)^{1/2} e^{-2kz\sqrt{2\theta}}}{4k^2-1} \,.
\]
All the terms in the sum are Laplace transforms and \cite{MR0352889}, p. 258, formula 5.87,
give the the inverses.  Since all the terms are positive the series can be inverted termwise
and we get
\begin{equation}
\label{FinalInversion}
\frac{1}{\sqrt{s}} P( M^+ - M^- \ge \frac{z}{\sqrt{s}}) =
\frac{1}{\sqrt{s}} \sum_{k=1}^\infty \frac{1}{4k^2-1} \cdot 
e^{-2k^2z^2/s} 
\end{equation}
or
\[
P(M^+ - M^- \ge z) = \sum_{k=1}^\infty \frac{1}{4k^2-1} \cdot
e^{-2k^2z^2} \,.
\]
Note that the distribution of $M^{+} - M^{-}$ is symmetric about $0$.


We now turn to the quotient  $Q=M^+/M^-$.  We can
multiply the numerator and denominator by $\sqrt{\gamma}$ and choose $\theta=1/2$.  Conditionally
on $L(S_\theta)=t$ an elementary calculation gives the conditional distribution function of $Q$
as
\[
F_{Q\vert L(S_\theta)=t}(z)
=
-\int_0^\infty \, tm^\prime(x) e^{-tm(x)} e^{-tm(xz)} \, \mathrm{d}x \,.
\]
Unconditioning and changing the order of integration gives
\begin{equation}
\label{IntegralDensityQ}
F_Q(z)
=
-\int_0^\infty \, \frac{m^\prime(x) \, \mathrm{d}x}{\left ( 1+m(x)+m(xz)\right )^2} \,.
\end{equation}
Substituting (\ref{Functionm}) and observing that
\[
m^\prime(x)
=
-\frac{1}{2 \, \mathrm{sinh}^2 x}
\]
gives
\begin{eqnarray*}
F_Q(z)\
& = &
\frac{1}{2} \, \int_0^\infty \, 
\frac{(1-e^{-2x})^2(1-e^{-2xz})^2 \, \mathrm{d}x}{\mathrm{sinh}^2 x 
\left ( 1-e^{-2x(z+1)} \right )^2 } \\
& = &
2 \, \int_0^\infty \, \frac{\mathrm{sinh}^2 zx \, \mathrm{d}x}{\mathrm{sinh}^2 \left ( (z+1)x \right ) } \\
& = &
\frac{1}{z+1} \, \left ( 1-\frac{\pi z \mathrm{cot}\left ( \frac{\pi z}{z+1} \right )}{z+1} \right ) \\
\end{eqnarray*}
where the last integral is given in \cite{MR1243179}, formula 3.511.9.  
\end{proof}




\bigskip

\leftline{\sc 3.3 Covariance and correlation.}

\medskip
The correlation between $M^+$ and
$M^-$ may be a quantity of interest.  Let
$U=\max_{0 \le t \le g_{S_\theta}} B_t$
and $V=-\min_{0 \le t \le g_{S_\theta}} B_t$.  The covariance of $U$ and
$V$ will be computed first.  By symmetry 
\[
E(U \vert L(S_\theta))=E(V \vert L(S_\theta)) \,.
\]
We have
\begin{eqnarray} \nonumber
\label{ConditionalCovariance}
\lefteqn{\qquad\mathrm{cov}(U,V)=} \\
& \phantom{=} &
\qquad = E \left (\mathrm{cov}(U,V \vert L(S_\theta) \right ) + 
\mathrm{cov} \left (E(U \vert L(S_\theta)), E(V \vert L(S_\theta) \right ) \,. \\  \nonumber
\end{eqnarray}
Denote $E(U \vert L(S_\theta))=\psi(L(S_\theta))$.  
By Theorem \ref{Triple} given $L(S_\theta)$ the conditional covariance of $U$ and $V$ 
is $0$.  It follows
\[ 
\mathrm{cov}(U,V)
=
\mathrm{var} \left ( \psi(L(S_\theta)) \right ) \,.
\]
The conditional expectation is computed from Theorem \ref{Triple}, (iii), as
\begin{equation}
\label{ConditionalExpectation}
\psi(t)=
\int_0^\infty (1-e^{-t m(x)}) \,\mathrm{d}x \\
\end{equation}
We can choose $2\theta=1$ so that the local time $L(S_\theta)$ is exponential with
parameter 1 and compute
\begin{eqnarray}
\label{ExpectationOfU}
E(U)
& = &
\int_0^\infty \, e^{-t} \psi(t) \, \mathrm{d}t \\ \nonumber
& = &
\int_0^\infty \, e^{-t} \, \mathrm{d}t \, \int_0^\infty (1-e^{-t m(x)}) \,\mathrm{d}x \\ \nonumber
& = &
\int_0^\infty \, \frac{m(x)}{1+m(x)} \,\mathrm{d}x \\ \nonumber
& = &
\int_0^\infty \, e^{-2x} \, \mathrm{d}x \\ \nonumber
& = &
1/2  \\ \nonumber
\end{eqnarray}
and
\begin{eqnarray}
\label{ExpectationOfPsiSquared} \nonumber
\lefteqn{E(\psi^2(L(S_{1/2})))=} \\ \nonumber
& = &
\int_0^\infty e^{-t} \,\mathrm{d}t
\int_0^\infty (1-e^{-t m(x)}) \,\mathrm{d}x
\int_0^\infty (1-e^{-t m(y)}) \,\mathrm{d}y \\ \nonumber
& = &
\int_0^\infty \,\mathrm{d}x \int_0^\infty \, \mathrm{d}y
\left [
1-\frac{1}{1+m(x)} -\frac{1}{1+m(y)} + \frac{1}{1+m(x)+m(y)}
\right ] \\ \nonumber
& = &
\int_0^\infty \,{\rm d}x \int_0^\infty \,{\rm d}y \;
\frac{e^{-2(x+y)} (2-e^{-2x}-e^{-2y})}{1-e^{-2(x+y)}} \\ \nonumber
& = &
\int_0^\infty \,{\rm d}t \frac{e^{-2t}}{1-e^{-2t}} 
\int_0^t (2-e^{-2(t-v)}-e^{-2v}) \,{\rm d}v \\ \nonumber
& = &
\int_0^\infty \frac{e^{-2t}}{1-e^{-2t}} (2t - (1-e^{-2t})) \,{\rm d}t \\ \nonumber
& = &
2 \int_0^\infty \sum_{k=1}^\infty t e^{-2kt} \,{\rm d}t - \frac{1}{2} \\ \nonumber
& = &
2\sum_{k=1}^\infty \frac{1}{4k^2} - \frac{1}{2} \\ \nonumber
& = &
\frac{\pi^2}{12} - \frac{1}{2} \\
\end{eqnarray}
The second line follows from Fubini's theorem and the third by changing
variables to $x+y=t$, $y=v$. It follows
\[
{\rm cov}(U,V) =
{\rm var}(\psi(L(S_\theta))) =
\frac{\pi^2}{12}-\frac{3}{4} \,.
\]
From (\ref{onesided}) one derives $E(M^+)=E(M^-) = \sqrt{2\pi}/4$ and
$\mathrm{var}(M^+)=\mathrm{var}(M^-)=1/2-\pi/8$.  By (\ref{BasicEqualityinLaw}) and
the independence of $\gamma$ and $(M^+,M^-)$ we have
\[
{\rm cov}(U,V) = E(\gamma) \, E(M^+ M^-) - E^2(\sqrt{\gamma}) \, E(M^+)E(M^-)
\]
we have
\[
E(M^+ M^-) = \frac{\pi^2}{12} - \frac{1}{2} \,.
\]
and
\begin{equation}
\label{CovarianceM+M-}
{\rm cov}(M^+,M^-) = \frac{\pi^2}{12} - \frac{1}{2} - \frac{\pi}{8} 
\end{equation}
yielding 
\begin{equation}
\label{CorrelationM+M-}
{\rm corr}(M^+,M^-) = 
\frac{\frac{\pi^2}{12} - \frac{1}{2} - \frac{\pi}{8}}{\frac{1}{2}-\frac{\pi}{8}} =
-0.654534 \,.
\end{equation}
\medskip

\par\noindent
{\bf Acknowledgements:}  
We owe thanks to a referee for several suggestions which 
improved the presentation.

\end{document}